\documentclass[12pt]{amsart} 
\usepackage{amsmath,amssymb}      
\usepackage{ascmac}
\usepackage{mathrsfs}
\usepackage[abbrev]{amsrefs}
\usepackage{comment}
\allowdisplaybreaks[1]
\newtheorem{thm}{Theorem}[section]

\newtheorem{prop}[thm]{Proposition}
\newtheorem{lem}[thm]{Lemma}

\newenvironment{pr}
   {{\noindent \bf Proof.   }}{\hfill \qed}

\def\<{\langle }

\newcommand{\eqsp}[1]{{\begin{equation}\begin{split}#1\end{split}\end{
equation}}}

\usepackage{graphicx,xcolor}

\begin{document}


\begin{center}
\textbf{\Large{$L^p$ estimate of the heat equation on a bounded domain}}
\end{center}
\vskip10mm

\centerline{Yoshinori Furuto${}^*$, Tsukasa Iwabuchi${}^{**}$ and Ryusei Kohama}
\vskip5mm
\centerline{Mathematical Institute, Tohoku University}
\centerline{Sendai 980-8578 Japan}
\footnote[0]{E-mail: $^*$yoshinori.furuto.p3@dc.tohoku.ac.jp,
$^{**}$t-iwabuchi@tohoku.ac.jp}

\begin{center}
 \begin{minipage}{120mm}
\small \textbf{Abstract.} 
We consider the linear heat equation on a bounded domain. 
We study estimates of the derivatives, up to the second order, of the solution locally in time in the Lebesgue spaces.
We give a self-contained proof of the estimates in the end-point cases $p = 1, \infty$.
 \end{minipage}
\end{center}

\section{Introduction}

Let $d \geq 2$. Suppose that $\Omega$ be a bounded domain of ${\mathbb R}^d$.
We consider the heat equation with the Dirichlet boundary condition.
\begin{equation} \label{main-eq}
    \left \{
        \begin{aligned}
            \partial_t u - \Delta u &= 0
                && \quad \text{ in } \ (0, T) \times \Omega, \\
            u &= 0
                && \quad \text{ on } \ (0, T) \times \partial \Omega, \\
            u(0, \cdot) &= u_0(\cdot)
                && \quad \text{ in } \ \Omega.
        \end{aligned}
    \right.
\end{equation}
Our aim is to obtain $L^p$ estimates of second-order derivatives.

Let us recall some classical results.
If $\Omega$ is a sufficiently smooth bounded domain and $p \in (1, \infty)$, then there exists the heat semigroup $\{S_p(t)\}_{t \geq 0}$ on $L^p(\Omega)$ and the solution $u(t) := S_p(t) u_0$, defined by the semigroup, satisfies the following $L^p$-$L^p$ estimates.
\begin{subequations} \label{20240409-1}
\begin{align}
    \| u(t) \|_{L^p} &\leq C \| u_0 \|_{L^p}, \label{20240409-1-a}
    \\
    \| \nabla u(t) \|_{L^p} &\leq C t^{-1/2} \| u_0 \|_{L^p}, \label{20240409-1-b}
    \qquad \text{ for } t > 0.
    \\
    \| \nabla^2 u(t) \|_{L^p} &\leq C t^{-1} \| u_0 \|_{L^p}, \label{20240409-1-c}
\end{align}
\end{subequations}
where $\| \cdot \|_{L^p}$ is the norm of $L^p(\Omega)$.
We notice that $u(t)$ decays exponentially in $L^p(\Omega)$ for large $t$. 
Our aim is to consider it locally in time.

In what follows, we suppose $u_0 \in C_c ^\infty (\Omega)$, which implies 
$S_2 (t) u_0 = S_p(t) u_0$ for any $p \in (1,\infty)$. 
We then abbreviate $p$, 
and often use the following notation  
\[
S(t)u_0 = S_p(t)u_0. 
\]

\begin{thm} \label{thm:p}
    Let $d \geq 2$. Suppose that $\Omega$ be a bounded domain of ${\mathbb R}^d$ with $C^3$ boundary, 
    $u_0 \in C_c^\infty (\Omega)$, 
    and $1 \leq p \leq \infty$. 
    Then $u(t) := S(t) u_0$ satisfies that
    \begin{equation} \label{eq:20240405-1}
    \begin{aligned} 
        \|  u(t) \|_{L^p} &\leq C \| u_0 \|_{L^p}, 
        \\
        \| \nabla u(t) \|_{L^p} &\leq C t^{-1/2} \| u_0 \|_{L^p}, 
        \\
        \| \nabla^2 u(t) \|_{L^p} &\leq C t^{-1} \| u_0 \|_{L^p},
    \end{aligned}
    \end{equation}
    for all $t > 0$.
\end{thm}

\noindent 
{\bf Remark. } 
The smoothness of the boundary of $\Omega$ is needed for the $L^\infty$ estimate, where we follow the argument of the paper by Abe-Giga~\cite{paper:AbGi-2013}.
We do not impose any assumption for the $L^1$ estimate. 
We will prove it by modifying the argument in~\cite{paper:IwMaTa-2018}.

We refer to some existing literature on derivative estimates.
In the case when $\Omega$ is a smooth bounded domain, the estimate of the solution itself is obvious from the property of the semigroup.
For the first derivatives, \eqref{20240409-1-b} holds for $p \in [1, \infty]$ (see~\cite{paper:Iw-2018},~\cite{paper:Shen-2005} and~\cite{paper:AbGi-2013}).
For the second derivatives, if $p \in (1, \infty)$, then~\eqref{20240409-1-c} can be obtained by using $L^p$ elliptic estimate and basic property of analytic semigroup (see~\cite{book:Lunardi} and~\cite{book:Gilbarg-Trudinger}). In addition,~\eqref{20240409-1-c} holds for $p = \infty$ (see~\cite{paper:AbGi-2013}).

We give a few comments on our proof.
We prove only the case $p = 1, \infty$.
The boundedness of without derivatives is due to the basic property of the $C_0$ semigroup.
For $p = 1$, we argue similarly to the paper~\cite{paper:IwMaTa-2018} who considered the first derivative. As for $p = \infty$, we apply the contradiction argument by~\cite{paper:AbGi-2013} for the Stokes equations.
Although not used in this paper, there is a possibility to show this theorem by direct estimate of the heat kernel using the method in \cite{paper:IsTa-2022}.

\vskip3mm 

This paper is organized as follows.
In Section 2, we prepare a lemma and a proposition to prove the theorem for $p = \infty$.
In Sections 3 and 4, we prove the theorem for $p = 1, \infty$.
In the Appendix, we show some lemmas for the paper to be self-contained.

\vskip3mm

\noindent {\bf Notation.} 
We denote by $A$ the Dirichlet Laplacian, which is defined by 
\[
\begin{cases}
D(A) = \{ f \in H^1_0 (\Omega) \, | \, \Delta f \in L^2 (\Omega)\}, 
\\
\displaystyle Af = - \sum _{ j=1}^d \partial _{x_j}^2 f , \quad f \in D(A) . 
\end{cases}
\]
We notice that $A$ is a self-adjoint operator on $L^2(\Omega)$ and by the argument in the paper ~\cite{paper:Iw-2018}, 
$e^{-tA}$ can be defined in more generalized distribution spaces 
based on the spectral multiplier. 
If $u_0 \in C_c^\infty (\Omega)$, then $S(t) u_0 = e^{-tA}u_0$, which is used in the proof for $p = 1$.

\section{Preliminaries}

We prepare a lemma and a proposition for the proof when $p = \infty$.
We consider the estimate for the following quantity.
\[
N[u](t, x) := |u(t, x)| + t^{1/2} |\nabla u(t, x)| + t |\nabla^2 u(t, x)| + t |\partial_t u(t, x)|.
\]

The first lemma considers an estimate of $N[u]$, and that for the Stokes equations is known~(see~\cite{paper:AbGi-2013}).
\begin{lem} \label{C2}
    Suppose that $\Omega$ is a bounded domain with $C^3$ boundary. Let $d < r < \infty$. Assume that $u_0 \in D(A_r)$, where $-A_r$ is the generator of the heat semigroup in $L^r(\Omega)$, and define $u(t) := S_r(t) u_0$.
    Then, $u(t, \cdot) \in C^2(\Omega)$ for all $t > 0$ and
    \[
    \sup_{0 < t < 1}{\| N[u](t) \|_{L^\infty} }< \infty.
    \]
\end{lem}

\vspace{3mm}

\begin{pr}
    We will show that there exists $C = C(\Omega, r) > 0$ such that
    \begin{equation} \label{LemmaGoal}
    \begin{split}
        \sup_{0 < t < 1}{(\| u(t) \|_{W^{1, r}} + t^{1/2} \| \nabla u(t) \|_{W^{1, r}} +  
        t \| \nabla^2 u(t) \|_{W^{1, r}} + t \| \partial_t u(t) \|_{W^{1, r}} )}
        \\
        \leq C \| u_0 \|_{D(A_r)},
    \end{split}
    \end{equation}
    where $\| u_0 \|_{D(A_r)} := \| u_0 \|_{L^r} + \| A_r u_0 \|_{L^r}$.

    By the property of the analytic semigroup $S_r(t) = e^{-t A_r}$, we have
    \[
    \begin{aligned}
        \sup_{0 < t < 1}{\| u(t) \|_{D(A_r)}} + \sup_{0 < t < 1}{t \| \partial_t u(t) \|_{D(A_r)}} \leq C \| u_0 \|_{D(A_r)},
    \end{aligned}
    \]
    where $C = C(\Omega, r)$. Thus we have proved
    \[
    \sup_{0 < t < 1}{(\| u(t) \|_{W^{1, r}} + \| \nabla u(t) \|_{W^{1, r}} + 
    t \| \partial_t u(t) \|_{W^{2, r}} )}
    \leq C \| u_0 \|_{D(A_r)}
    \]
    as the $D(A_r)$ norm and the $W^{2, r}(\Omega)$ norm are equivalent by the elliptic estimate.
    The terms
    \[
    \| u(t) \|_{W^{1, r}}, \quad
    t^{1/2} \| \nabla u(t) \|_{W^{1, r}}, \quad 
    t \| \nabla^2 u(t) \|_{L^r}, \quad
    t \| \partial_t u(t) \|_{W^{1, r}}
    \]
    are less than or equal to the left hand side above.


    It remains to prove that
    \[
    \sup_{0 < t < 1}{t \| \nabla^3 u(t) \|_{L^r}}
    \leq C \| u_0 \|_{D(A_r)}.
    \]
    Applying the estimate in a nonhomogeneous elliptic equation (Theorem~\ref{high-est-thm} for $k = 1$ and $f = \partial_t u$),
    we can obtain
    \[
    \begin{aligned}
        \| \nabla^3 u(t) \|_{L^r}
        &\leq C ( \| u(t) \|_{L^r} + \| \partial_t u(t) \|_{W^{1, r}}) 
        \\
        &\leq C \| u_0 \|_{D(A_r)}, 
        \quad t > 0, 
    \end{aligned}
    \]
    which proves the above inequality.
    Here the assumption that $\Omega$ is a $C^3$ domain is needed for Theorem~\ref{high-est-thm}.
\end{pr}

\begin{prop}\label{prop:0126-1}
Let $\Omega = \mathbb R^d_+$. 
Suppose that $u$ satisfies
\[
\int_0^T \int_{\Omega} u (\partial_t \phi + \Delta \phi) dx dt 
= 0
\]
for any
$
\phi \in C([0, T] ; W^{1,1}_0(\Omega) \cap W^{2, 1}(\Omega) \cap C(\Omega)) 
\cap C^1([0, T] ; L^1(\Omega))
$
with $| \nabla \phi(x)| \leq C(1 + |x|^2)^{-d/2}$.
Then $u$ is $0$ almost everywhere.
\\
\end{prop}

\section{$L^1$ estimate} 

We follow the argument as in the paper~\cite{paper:IwMaTa-2018}. 
Instead of $S(t)u_0$, we consider more general problem, which is the boundedness of the 
spectral multiplier $\varphi(t A)$, 
where $\varphi : \mathbb R \to \mathbb R$ 
with ${\rm supp \, }\varphi \subset [0,\infty)$. 

Following the proof of Theorem~2.1 in \cite{paper:IwMaTa-2018}, 
we decompose 
$\Omega$ by using cubes whose side length 
is $t^{\frac{1}{2}}$, and apply the H\"older 
inequality to have that 
\begin{equation}\label{eq:0131-1}
\| \nabla ^2 \varphi (tA) f \|_{L^1} 
\leq t^{\frac{d}{4}} 
\| \nabla ^2 \varphi(tA) f \|_{\ell^1 (L^2)_t},
\end{equation}
where 
\[
\| \nabla ^2 \varphi(tA) f \|_{\ell^1 (L^2)_t}
:= \sum_{n \in {\mathbb Z}^d} \| \nabla ^2 \varphi(tA) f \|_{L^2 (C_t(n))}, 
\]
$C_t(n)$ is a cube with side $t^{1/2}$ centered $t^{1/2} n \in {\mathbb R}^d$.
Introduce $\widetilde \varphi \in \mathcal S(\mathbb R)$ such that  
\[
\widetilde \varphi(\lambda) 
= (1 + \lambda) ^\beta \varphi(\lambda), 
\quad \beta > 0 . 
\]
We can write 
\begin{equation}\notag 
\begin{split}
& 
\| \nabla ^2 \varphi(tA) f \|_{\ell^1 (L^2)_t}
\\
=& \| \nabla ^2 \widetilde \varphi(tA) 
(1+tA)^{-\beta} f \|_{\ell^1 (L^2)_t}
\\
\leq & C 
\Big(  \| \nabla ^2 \widetilde \varphi(tA) \| _{L^2 \to L^2} 
+ t ^{-\frac{d}{4}} 
 ||| \nabla ^2 \widetilde \varphi(tA) |||_\alpha 
  ^{\frac{d}{2\alpha}} 
 \| \nabla ^2 \widetilde \varphi(tA) \| _{L^2 \to L^2} ^{1-\frac{d}{2\alpha}} 
\Big)
\| (1+tA)^{-\beta} f \|_{\ell^1 (L^2)_t}
\end{split}
\end{equation}
for any $\alpha > d/2$, where 
\[
||| \nabla ^2 \widetilde \varphi(tA) |||_\alpha :=
\sup_{n \in {\mathbb Z}^d} \left\| | \cdot - t^{1/2} n |^\alpha \nabla ^2 \widetilde \varphi(tA) \chi_{C_t(n)}  \right\|_{L^2 \to L^2} . 
\]
The elliptic estimate in $L^2 (\Omega)$ shows that 
\[
\| \nabla ^2 \widetilde \varphi(tA) \| _{L^2 \to L^2} 
\leq C \| A \widetilde \varphi(tA) \| _{L^2 \to L^2} 
\leq C t^{-1}. 
\]
As for 
$ ||| \nabla ^2 \widetilde \varphi(tA) |||_\alpha $, 
by following the proof of Lemma~7.1 
in \cite{paper:IwMaTa-2018}, 
$\nabla ^2 \widetilde \varphi(tA)$ can be 
written with the following 
\[
\begin{aligned}
    \nabla ^2 (1+tA)^{-1}, \quad 
    \nabla (1+tA)^{-1}, \quad (1+tA)^{-1}, \quad 
    e^{is (1+tA)^{-1}},
    \\
    \text{integrals of Fourier transform of } \widetilde \varphi,
\end{aligned}
\]
and we conclude that 
\[
||| \nabla ^2 \widetilde \varphi(tA) |||_\alpha 
\leq C t ^{\frac{\alpha-2}{2}} . 
\]
The above two inequlities show that 
\begin{equation}\label{eq:0131-2}
\begin{split}
\| \nabla ^2 \varphi(tA) f \|_{\ell^1 (L^2)_t}
\leq 
& 
C(t^{-1} + t^{-\frac{d}{4}} (t ^{\frac{\alpha-2}{2}})^{\frac{d}{2\alpha}}
(t^{-1})^{1-\frac{d}{2\alpha}}) 
\| (1+tA)^{-\beta} f \|_{\ell^1 (L^2)_t} 
\\
= & C t^{-1} 
\| (1+tA)^{-\beta} f \|_{\ell^1 (L^2)_t} .
\end{split}
\end{equation}
Finally, it follows by 
Proposition~7.1 in \cite{paper:IwMaTa-2018} 
that 
\[
\| (1+tA)^{-\beta} f \|_{\ell^1 (L^2)_t}
\leq C t^{-\frac{d}{4}} \| f \|_{L^1}.  
\]
We then obtain by \eqref{eq:0131-1} that 
\[
t^{\frac{d}{4}} \| \nabla ^2 \varphi(tA) f \|_{\ell^1 (L^2)_t}
\leq C t^{\frac{d}{4}} t^{-1} \| (1+tA)^{-\beta} f \|_{\ell^1 (L^2)_t}
\leq C t^{-1} \| f \|_{L^1}. 
\]

\section{$L^\infty$ estimate}

\begin{thm}
    Let $\Omega \subset {\mathbb R}^d$ be a bounded domain with $C^3$ boundary.
    Then, there exists $C > 0$ and $T > 0$ such that
    \[
    \sup_{0 < t < T}{\| N[u](t) \|_{L^\infty}} < C \| u_0 \|_{L^\infty}
    \]
    holds for any solution $u(t)$ of $\eqref{main-eq}$ with $u_0 \in C_c^{\infty}(\Omega).$
\end{thm}

\begin{pr}
    Prove by contradiction. Assume that there exist a sequence $\{ u_m \}_{m \in {\mathbb N}}$ of solution to $\eqref{main-eq}$ and a sequence $\{ \tau_m \}_{m \in {\mathbb N}}$ such that
    \[
    \begin{aligned}
        \| N[u](\tau_m) \|_{L^\infty} 
        > m \| u_{0, m} \|_{L^\infty}, 
        \\
        \tau_m \searrow 0 \, (m \to \infty)
    \end{aligned}
    \]
    where $u_{0, m} \in C_c^{\infty}(\Omega)$ is an initial data of $u_m$.

    Now we introduce a function $v_m$ as follows.
    It holds that
    $\displaystyle M_m := \sup_{t \in (0, \tau_m)}{\| N[u_m](t) \|_{L^\infty}} < \infty$ 
    by Lemma $\ref{C2}$, and we define $\tilde{u}_m := u_m / M_m$.
    Then, there exists some point $(t_m, x_m) \in (0, \tau_m) \times \Omega$ satisfying
    $N[\tilde{u}_m](t_m, x_m) \geq 1/2$.
    We therefore define $v_m$ by translating and scaling where the point $(t_m, x_m)$ corresponds to the point $(1, 0)$.
    \begin{align*}
        v_m (t, x) := \tilde{u}_m (t_m t, x_m + t_m^{1/2} x) 
            \quad \text{ for } (t, x) \in (0, 1) \times \Omega_m, \\
        \text{where } \Omega_m := \Big\{ y \in \mathbb{R}^d \, \Big| \, y = \frac{x - x_m}{t_m^{1/2}}, x \in \Omega \Big\}.
    \end{align*}

    Here $v_m$ satisfies the heat equation and the following inequities.
    \[
    \begin{aligned}
        \mathrm{N}[v_m](1, 0) = 
            \mathrm{N}[\tilde{u}_m](t_m, x_m) \geq \frac{1}{2}, 
        \\
        \sup_{t \in (0, 1)} \| \mathrm{N}[v_m](t) \|_{L^\infty} = 
            \sup_{t \in (0, t_m)} \| \mathrm{N}[\tilde{u}_m](t) \|_{L^\infty} \leq 1, 
        \\
        \| v_{0, m} \|_{L^\infty}  = 
            \| \tilde{u}_{0, m} \|_{L^\infty} < \frac{1}{m}.
    \end{aligned}
    \]

    The situation depends on whether the limit of
    \[
    c_m := \frac{d(\Omega, x_m)}{t_m^{1/2}} = d(\Omega_m, 0)
    \]
    goes infinity or not.

    \vspace{3mm}
    
    Case 1. $\displaystyle \limsup_{m \to \infty}c_m = \infty$

    Take the subsequence so that $\displaystyle \lim_{m \to \infty} c_m = \infty$.
    We denote $\tilde{v}_m$ as an extension of $v_m$ to ${\mathbb R}^d$. 

    By a compact embedding $W^{3, r}(\Omega) \subset \subset W^{2, \infty}(\Omega)$, 
    there exist a subsequence $\{ \tilde{v}_{m(k)} \}_{k \in \mathbb N}$ and $\tilde{v}$ such that 
    $\tilde{v}_{m(k)}$ converges to $\tilde{v}$ locally uniformly in ${\mathbb R}^d \times (0, 1]$.

    For any $\phi \in C_c^{\infty} ({\mathbb R}^d \times [0, 1))$, it follows
    \[
        \int_0^1 \int_{{\mathbb R}^d} \tilde{v}_{m(k)} \cdot (\partial_t \phi + \Delta \phi) dx dt
        = \int_{{\mathbb R}^d} \tilde{v}_{m(k)} (x, 0) \cdot \phi(x, 0) dx, 
    \]
    and by taking the limit as $k \to \infty$, we obtain
    \[
    \int_0^1 \int_{{\mathbb R}^d} \tilde{v} \cdot (\partial_t \phi + \Delta \phi) dx dt
    = 0.
    \]
    Therefore the limit $\tilde{v}$ is the weak solution of heat equation and 
    it follows $\tilde{v} \equiv 0$ by the uniqueness of the heat equation (see section~4.4.2 in~\cite{book:Giga-Giga}).
    However, it contradicts $N[\tilde{v}](1, 0) \geq 1/2$.

    \vspace{3mm}

    Case 2. $\displaystyle \limsup_{m \to \infty}c_m < \infty$

    By taking a subsequence, we assume that $\displaystyle \lim_{m \to \infty} c_m = c_0 < \infty$ and 
    $x_m$ converges to $\hat{x} \in \partial\Omega$.
    In addition, by rotation and translation, we also assume $\hat{x} = 0$.

    Since $\Omega$ is of $C^3$, there exisits a $C^3$ function $h$ such that 
    \[
    \begin{aligned}
        (\Omega_m)_{\text{loc}} := 
        \{ (y', y_d) \in \mathbb{R}^d \mid
        &h(x_m' + t_m^{1/2} y') < (x_m)_d + t_m^{1/2} y_d
        \\
        &< h(x_m' + t_m^{1/2} y') + \beta, 
        |t_m^{1/2} y'| < \alpha \}
        \subset \Omega_m, 
        \\
        \nabla{h}(0) = 0, h(0) = 0. 
    \end{aligned}
    \]
    where $\alpha, \beta$ are the constants independent of $\hat{x} \in \partial\Omega$.

    Since $d(\Omega_m, x_m) / (x_m)_d \to 1$ as $m \to \infty$, 
    the domain $(\Omega_m)_{\text{loc}}$ approaches to 
    \[
    \mathbb{R}^d_{+, -c_0} := \{ (x', x_d) \in \mathbb{R}^d \mid x_d > -c_0 \}.
    \]


    Fix $R > 0$ and let $B_R^+ = B_R(0, \ldots, 0, -c_0) \cap {\mathbb R}^d_{+, -c_0}$.

    Define $\tilde{v}_m \in W^{3, r}(\mathbb R^d \times (0, 1])$ as an extension of $v_m$ and there exists $\tilde{v} \in W^{2, \infty}({\mathbb R}^d_{+, -c_0} \times (0, 1])$ and $C > 0$ such that
    \[
    \begin{aligned}
    \tilde{v}_{m(k)} \to \tilde{v}&
    \\ 
    \nabla \tilde{v}_{m(k)} \to \nabla \tilde{v}&
    \quad \text{ locally uniformly in } \overline{\mathbb{R}}^d_{+, -c_0} \times (0, 1], 
    \\
    \nabla^2 \tilde{v}_{m(k)} \to \nabla^2 \tilde{v}&
    \end{aligned}
    \]
    \[
    \| \tilde{v}_m \|_{W^{2, \infty}({\mathbb R^d_{+, -c_0}} \times (0, 1))} \leq 
    C \| v_m \|_{W^{2, \infty}({B_R^+} \times (0, 1))}
    \text{ for all } m \in \mathbb N
    \]
    where $\tilde{v}_{m(k)}$ is a subsequence of $\{ \tilde{v}_m \}_{m \in \mathbb N}$.
    For simplicity, we denote the subsequence by $\{ \tilde{v}_m \}_{m \in \mathbb N}$.

    Fix a test funciton 
    \[
    \phi \in  C\Big([0, 1] ; W^{1,1}_0(R^d_{+, -c_0}) \cap W^{2, 1}(R^d_{+, -c_0}) \cap C(R^d_{+, -c_0})\Big) \cap C^1\Big([0, 1] ; L^1(R^d_{+, -c_0})\Big)
    \]
    with $|\nabla \phi(x)| \leq C (1 + |x|^2)^{-d/2}$.

    
    Since $v_m$ satisfies $\eqref{main-eq}$ in $\Omega_m \times (0, 1]$, it follows
    \[
    \begin{aligned}
    0 
    &= \int_0^1 \int_{\Omega_m \cap B_R^+} \Delta \tilde{v}_m \phi dx dt
    - \int_0^1 \int_{\Omega_m \cap B_R^+} (\partial_t \tilde{v}_m) \phi dx dt
    \\
    &= - \int_0^1 \int_{\Omega_m \cap B_R^+} \nabla \tilde{v}_m \cdot \nabla \phi dx dt
    + \int_0^1 \int_{\partial(\Omega_m \cap B_R^+)} (\nabla \tilde{v}_m \cdot \boldsymbol{n}) \phi dx dt
    \\
    &\, - \int_0^1 \int_{\Omega_m \cap B_R^+} *\partial_t \tilde{v}_m) \phi dx dt.
    \end{aligned}
    \]
    For the second term, divide the integral domain into $\partial(\Omega_m \cap B_R^+) \cap \partial B_R^+$ and $\partial(\Omega_m \cap B_R^+) \setminus \partial B_R^+$. Then, by the decay of  $| \phi |$, the continuity of $\phi$ and $\phi \in W_0^{1, 1}$, it holds
    \[
    \begin{aligned}
    \sup_{m \in \mathbb N} \left| \int_0^1 \int_{\partial(\Omega_m \cap B_R^+) \cap \partial B_R^+} (\nabla \tilde{v}_m \cdot \boldsymbol{n}) \phi dx dt \right| \to 0
    \quad \text{ as } R \to \infty, 
    \\
    \left| \int_0^1 \int_{\partial(\Omega_m \cap B_R^+) \setminus \partial B_R^+} (\nabla \tilde{v}_m \cdot \boldsymbol{n}) \phi dx dt \right| \to 0
    \quad \text{ as } m \to \infty.
    \end{aligned}
    \]
    From the Lebesgue's dominated convergence theorem and Fubini's theorem the integration respect with t, for any $\delta > 0$, there exists $M \in \mathbb N$ such that for $m \geq M$,
    \[
    \begin{aligned}
        & \left| \int_0^1 \int_{\Omega_m \cap B_R^+} \nabla \tilde{v}_m \cdot \nabla \phi dx dt
        - \int_0^1 \int_{B_R^+} \nabla \tilde{v}_m \cdot \nabla \phi dx dt \right| < \delta,
        \\
        & \left| \int_0^1 \int_{\Omega_m \cap B_R^+} (\partial_t \tilde{v}_m) \phi dx dt
        - \int_0^1 \int_{B_R^+} (\partial_t \tilde{v}_m) \phi dx dt \right| < \delta.
    \end{aligned}
    \]
    We consider the integral on $B_R^+$. By integration by parts, 
    \begin{equation}
    \begin{split} \notag
    &- \int_0^1 \int_{B_R^+} \nabla \tilde{v}_m \cdot \nabla \phi dx dt
    \\
    = &\int_0^1 \int_{B_R^+} \tilde{v}_m \Delta \phi dx dt
    - \int_0^1 \int_{\partial B_R^+} \tilde{v}_m (\nabla \phi \cdot \boldsymbol{n}) dx dt.
    \end{split}
    \end{equation}
    For the second term, we again divide the integral domain into $\partial(\Omega_m \cap B_R^+) \cap \partial B_R^+$ and $\partial(\Omega_m \cap B_R^+) \setminus \partial B_R^+$. Then the decay of $| \nabla \phi| $ implies
    \[
    \begin{aligned}
    \sup_{m \in \mathbb N} \left| \int_0^1 \int_{\{x_d \neq -c_0\} \cap \partial B_R^+} \tilde{v}_m (\nabla \phi \cdot \boldsymbol{n}) dx dt \right| \to 0
    \quad \text{ as } R \to \infty.
    \end{aligned}
    \]
    Now we estimate $\tilde{v}_m(t, x)$ at $x \in \partial B_R^+, x_d = -c_0$.
    Let $w_m := \displaystyle \sup_{y \in \partial(\Omega_m \cap B_R^+) \setminus \partial B_R^+} | y_d |$ and notice that
    \[
    \begin{aligned}
    w_m \to 0 
    &\text { as } m \to \infty, 
    \\
    dist(x, \Omega_m) < w_m 
    &\text{ for all } x \in \partial B_R^+, x_d = -c_0. 
    \end{aligned}
    \]
    For any points $x \in \partial B_R^+, x_d = -c_0$, define $x_b$ as the point in $\partial \Omega_m$ that is closest to $x$, then
    \[
    \begin{aligned}
    |\tilde{v}_m(t, x)| 
    &=
    |\tilde{v}_m(t, x) - \tilde{v}_m(t, x_b)|
    \\
    &\leq
    \int_0^1 | \partial_\theta \tilde{v}_m(t, \theta x + (1 - \theta) x_b) | \, d\theta
    \\
    &\leq
    | x - x_b | \times \left| D \tilde{v}_m(t, \theta_0 x + (1 - \theta_0) x_b \right|
    \quad \text{ for some } \theta_0 \in (0, 1)
    \\
    &\leq
    w_m \times C t^{-1/2}.
    \end{aligned}
    \]
    By this estimate and the limit $w_m \to 0$ as $m \to \infty$, we obtain
    \[
    \int_0^1 \int_{\{x_d = -c_0\} \cap \partial B_R^+} \tilde{v}_m (\nabla \phi \cdot \boldsymbol{n}) dx dt \to 0
    \quad \text{ as } m \to \infty.
    \]

    From above, we take the limits $m \to \infty, R \to \infty$ in this order, then
    \[
    \lim_{R \to \infty} \lim_{m \to \infty} \int_0^1 \int_{\Omega_m \cap B_R^+} \Delta \tilde{v}_m \phi dx dt
    = \int_0^1 \int_{\mathbb R^d_{+, -c_0}} \tilde{v} \Delta \phi dx dt.
    \]
    Similarly, it holds
    \[
    \begin{aligned}
    \lim_{R \to \infty} \lim_{m \to \infty} \int_0^1 \int_{\Omega_m \cap B_R^+} \partial_t \tilde{v}_m \phi dx dt
    = - \int_0^1 \int_{\mathbb R^d_{+, -c_0}} \tilde{v} \partial_t \phi dx dt
    - \int_{\mathbb R^d_{+, -c_0}} \tilde{v}(0, x) \phi(0, x) dx,
    \end{aligned}
    \]
    and therefore, $\tilde{v}$ satisfies
    \[
    \int_0^1 \int_{\mathbb R^d_{+, -c_0}} \tilde{v} (\Delta \phi +\partial_t \phi) dx dt
    = - \int_{\mathbb R^d_{+, -c_0}} \tilde{v}(0, x) \phi(0, x) dx.
    \]
    In other words, $\tilde{v}$ is a weak solution of the heat equation on ${\mathbb R}^d_{+, -c_0}$.

    By the uniqueness of the heat equation for the half space (Proposition~\ref{prop:0126-1}), it follows $\tilde{v} \equiv 0$.
    However, it contradicts $N[\tilde{v}](1, 0) \geq 1/2$.
\end{pr}

\appendix
\section{higher-order derivative estimates}

We show the following theorem for the paper to be self-contained. 
Proof idea can be found in the book by Gilbarg-Trudinger~\cite{book:Gilbarg-Trudinger} (Theorem 9.19.) .

\begin{thm} \label{high-est-thm}
    Suppose $k \in \mathbb N \cup \{0\}$. Let $\Omega \subset {\mathbb R}^d$ be a bounded domain with $C^{k+2}$ boundary and $p \in (1, \infty)$.
    If $u \in W^{2, p}(\Omega), f \in W^{k, p}(\Omega)$ satisfy
    \begin{equation}
        \left \{
            \begin{aligned}
                \Delta u &= f
                    && \quad \textup{ in } \Omega, \\
                u &= 0
                    && \quad \textup{ on } \partial \Omega,
            \end{aligned}
        \right.
    \end{equation}
    then $u \in W^{k+2, p}(\Omega)$ and there exists $C = C(n, k, p)$ such that
    \[
    \| D^{k+2} u \|_{L^p(\Omega)} \leq C ( \| u \|_{L^p(\Omega)} + \| f \|_{W^{k, p}(\Omega)} )
    \]
\end{thm}

The proof is derived from an internal estimate and a boundary estimate.

\begin{lem}
    Assume the same as in the Theorem $\ref{high-est-thm}$.
    Then, for any subdomain $\Omega' \subset\subset \Omega$, $u \in W^{k+2, p}(\Omega')$ and satisfies
    \[
    \| D^{k+2} u \|_{L^p(\Omega')} \leq C ( \| u \|_{L^p(\Omega)} + \| f \|_{W^{k, p}(\Omega)} ).
    \]
\end{lem}

\begin{pr}
    We refer to Theorem 9.11. in ~\cite{book:Gilbarg-Trudinger} in the case when $k = 0$. 

    Let $k = l \geq 1$. For any $i = 1, \ldots, d$, $u$ and $f$ satisfy
    \[
    \Delta (\partial_i u) = \partial_i f
    \quad \text{ in } \Omega, 
    \]
    and 
    \[
    \| D^{l+1} (\partial_i u) \|_{L^p(\Omega')}
    \leq C ( \| \partial_i u \|_{L^p(\Omega)} + \| \partial_i f \|_{W^{l-1, p}(\Omega)} ).
    \]
     by the estimate for $k = l - 1$.

     The elliptic estimiate implies
     \[
     \| D u \|_{L^p(\Omega)} \leq C \left( \| u \|_{L^p(\Omega)} + \| f \|_{L^p(\Omega)} \right)
     \]
    and therefore, 
    \[
    \begin{aligned}
        \| D^{l+2} u \|_{L^p(\Omega')}
        \leq& C ( \| Du \|_{L^p(\Omega)} + \| Df \|_{W^{l-1, p}(\Omega)} )
        \\
        \leq& C ( \| u \|_{L^p(\Omega)} + \| f \|_{W^{l, p}(\Omega)} ).
    \end{aligned}
    \]
\end{pr}

\begin{lem} \label{lem:0310-1}
    Let $\Omega^{+} := \Omega \cap \mathbb{R}^d_+, (\partial\Omega)^+ := \partial\Omega \cap \mathbb{R}^d_+$.
    If $u \in W^{1, p}_0(\Omega^+), f \in W^{k, p}(\Omega^+)$ satisfy
    \begin{equation}
        \left \{
            \begin{aligned}
                \Delta u &= f
                    && \quad \textup{ in } \Omega^+, \\
                u &= 0
                    && \quad \textup{ on } (\partial \Omega)^+,
            \end{aligned}
        \right.
    \end{equation}
    then $u \in W^{k+2, p}(\Omega^+)$ and there exists $C = C(n, k, p)$ such that
    \[
    \| D^{k+2} u \|_{L^p(\Omega^+)} \leq C \| D^k f \|_{L^p(\Omega^+)}.
    \]
\end{lem}

\begin{pr}
    We follow the argument as Lemma 9.12 in ~\cite{book:Gilbarg-Trudinger}. 

    Let $k \geq 1$. Extend $u$ and $f$ to $\mathbb R^d_+$ by zero-extension and then to $\mathbb R^d$ by odd reflection. Then $u \in W^{1, p}(\mathbb R^d)$ is a weak solution of $\Delta u = f$.

    Since $\text{supp } u \subset \mathbb R^d$ is compact, the regularization $u_h := u * \rho_h \in C^\infty_c(\mathbb R^d), f_h := f * \rho_h \in C^\infty_c(\mathbb R^d)$ where $\rho_h$ is a mollifier such that 
    \[
    \begin{aligned}
        &\text{ supp } \rho_h \subset \{ |x| < h \}, \\
        &\rho_h * g \to g \text{ in } L^p(\mathbb R^d) \quad (h \to 0)
    \end{aligned}
    \]
    for any $g \in L^p(\mathbb R^d)$.

    Then, it holds that $\Delta \partial_i u_h = \partial_i f_h$ and
    \[
    \begin{aligned}
    \| D^{k+2} u_h \|_{L^p(\mathbb R^d)}
    &\leq C \| \nabla^2 D^k u_h \|_{L^p(\mathbb R^d)}
    \\
    &\leq C \| \Delta D^k u_h \|_{L^p(\mathbb R^d)}
    \\
    &= C \| D^k f_h \|_{L^p(\mathbb R^d)}.
    \end{aligned}
    \]
    by Calderon-Zygmund inequality (see Theorem 9.4. in ~\cite{book:Gilbarg-Trudinger}). We take the limit as $h \to 0$, then 
    \[
    \| D^{k+2} u \|_{L^p(\mathbb R^d)}
    \leq C \| D^{k} f \|_{L^p(\mathbb R^d)}
    \]
    and
    \[
    \begin{aligned}
    \| D^{k+2} u \|_{L^p(\Omega)}
    &\leq \| D^{k+2} u \|_{L^p(\mathbb R^d)}
    \\
    &\leq C \| D^{k} f \|_{L^p(\mathbb R^d)}
    \\
    &\leq C \| D^{k} f \|_{L^p(\Omega)}
    \end{aligned}
    \]
    from the definition of the extension.
\end{pr}

\begin{lem}
    Assume the same as in the Theorem $\ref{high-est-thm}$.
    Then, for any subdomain $\Omega' \subset\subset \Omega \cup \partial \Omega$, $u \in W^{k+2, p}(\Omega')$ and satisfies
    \[
    \| D^{k+2} u \|_{L^p(\Omega')} \leq C ( \| u \|_{L^p(\Omega)} + \| f \|_{W^{k, p}(\Omega)} )
    \]
\end{lem}

\begin{pr}
    The argument is based on change of variables and Lemma \ref{lem:0310-1}, and same as Theorem 9.13 in the book~\cite{book:Gilbarg-Trudinger}. 
\end{pr}

\section{Proof of Proposition~\ref{prop:0126-1}}



\begin{pr}
The argument is based on the book by Giga~\cite{book:Giga-Giga}.

For every $f \in C_c^\infty ((0,T) \times \mathbb R^d_+)$, we consider
\begin{equation}
    \label{eq:0306-2} \notag
    \left \{
    \begin{aligned}
        \partial_t \psi - \Delta \psi &= f
            && \quad \text{ in } (0, T] \times \mathbb R^d_+
        \\
        \psi(0, x) &= 0
            && \quad \text{ on } \{t = 0\} \times \mathbb R^d_+.
    \end{aligned}
    \right.
\end{equation}

There exists a unique solution
\[
\psi(t) := \int_0^t e^{-(t-s)A} f(s) ds, 
\]
where
\[
\begin{aligned}
    e^{-(t-s)A} f = \int_{\mathbb R^d_+} (G_{t-s}(x'-y', x_d-y_d) - G_{t-s}(x'-y', x_d+y_d) ) f(s, y) dy, 
    \\
    G_t(x, y) = \frac{1}{(4 \pi t)^{d/2}} e^{- | x - y |^2 / 4t}.
\end{aligned}
\]

Since $f \in C_c^\infty ((0,T) \times \mathbb R^d_+)$, 
\[
\begin{aligned}
\psi = \psi_f \in C([0, T] ; W^{1,1}_0(R^d_+) \cap W^{2, 1}(R^d_+) \cap C(R^d_+)) 
\cap C^1([0, T] ; L^1(R^d_+)), 
\\
|\nabla \psi(x)| \leq C (1 + |x|^2)^{-d/2}
\end{aligned}
\]
(See~\cite{book:Cazenave-Haraux}).
We define $\widetilde \psi_f (t, x) := \psi(T-t, x)$ and then, for every $f \in C_c^\infty ((0,T) \times \mathbb R^d_+)$, 
\[
\begin{aligned}
    0
    &= \int_0^T \int_{\mathbb R^d_+} u (\partial_t \widetilde \psi_f + \Delta \widetilde \psi_f) dx dt
    \\
    &= \int_0^T \int_{\mathbb R^d_+} u (-\partial_t \psi_f + \Delta \psi_f) dx dt
    \\
    &= - \int_0^T \int_{\mathbb R^d_+} u f dx dt, 
\end{aligned}
\]
which implies $u \equiv 0$.
\end{pr}

\end{document}